\documentclass{amsart}
\usepackage{amsfonts}
\usepackage{amscd}
\usepackage{amssymb}
\usepackage{mathrsfs}
\theoremstyle{plain}
\newtheorem{thm}{Theorem}[section]
\newtheorem{cor}[thm]{Corollary}
\newtheorem{lem}[thm]{Lemma}
\newtheorem{pro}[thm]{Proposition}
\theoremstyle{definition}

\newtheorem{rmk}[thm]{Remark}
\newtheorem{defn}[thm]{Definition}
\usepackage{amsmath}
\usepackage{amssymb}
\usepackage{amsthm}
\usepackage{latexsym}
\usepackage[T1]{fontenc}
\usepackage[all]{xy}
\DeclareMathOperator{\id}{id}
\DeclareMathOperator{\Per}{Per} \DeclareMathOperator{\supp}{supp}
\DeclareMathOperator{\Ker}{Ker} 
 \DeclareMathOperator{\sep}{Sep}

\begin{document}
\author{Christian Svensson}
\address{Mathematical Institute, Leiden University,
P.O. Box 9512, 2300 RA Leiden, The Netherlands, and Centre for
Mathematical Sciences, Lund University, Box 118, SE-221 00 Lund,
Sweden} \email{chriss@math.leidenuniv.nl}
\author{Sergei Silvestrov}
\address{Centre for Mathematical Sciences,
Lund University, Box 118, SE-221 00 Lund, Sweden}
\email{Sergei.Silvestrov@math.lth.se}
\author{Marcel de Jeu}
\address{Mathematical Institute,
Leiden University, P.O. Box 9512, 2300 RA Leiden, The Netherlands}
\email{mdejeu@math.leidenuniv.nl}

\title[Dynamical systems associated with crossed products]{Dynamical
 systems associated with crossed products}


\noindent 
\noindent 
\\

\begin{abstract}
In this paper, we consider both algebraic crossed products of
 commutative complex algebras $A$ with
the integers under an automorphism of $A$, and Banach algebra crossed
 products of commutative $C^*$-algebras
$A$ with the integers under an automorphism of $A$. We investigate, in
 particular, connections between algebraic properties of
these crossed products and topological properties of naturally
 associated dynamical systems. For example, we draw conclusions about
the ideal structure of the crossed product by investigating the
 dynamics of such a system.
To begin with, we recall results in this direction in the context of
 an
 algebraic crossed product
 and give simplified proofs of generalizations of some of these
 results. We also investigate new questions, for
example about ideal intersection properties of algebras properly
 between the coefficient algebra $A$ and its commutant $A'$.
 Furthermore, we
introduce a Banach algebra crossed product and study the relation
 between the structure of this algebra and the topological dynamics of
 a naturally associated system.
\end{abstract}
\maketitle
\section{Introduction}
A lot of work has been done on the connection between
certain topological dynamical systems and crossed product
$C^*$-algebras. In \cite{tom1} and \cite{tom2}, for example, one
starts with a homeomorphism $\sigma$ of a compact Hausdorff space
$X$ and constructs the crossed product $C^*$-algebra $C(X)
\rtimes_{\alpha} \mathbb{Z}$, where $C(X)$ is the algebra of
continuous complex valued functions on $X$ and $\alpha$ is the
 $\mathbb{Z}$-action on $C(X)$ naturally induced by $\sigma$. One of
 many
results obtained is equivalence between simplicity of the algebra
and minimality of the system, provided that $X$ consists of
infinitely many points, see \cite{davidson}, \cite{power},
\cite{tom1}, \cite{tom2} or, for a more general approach in the
 metrizable case,
\cite{williams}. In \cite{svesildej}, a purely algebraic variant of
the crossed product is considered, having more general classes of
algebras than merely continuous functions on compact Hausdorff
spaces as coefficient algebras.
For example, it is proved there that, for such crossed products,
the analogue of the equivalence between density of non-periodic
points of a dynamical system and maximal commutativity of the
 coefficient algebra in the associated crossed product
$C^*$-algebra is true for significantly larger classes of
coefficient algebras and associated dynamical systems.
In \cite{art2}, further work is done in this setup, mainly for crossed
 products
of complex commutative semi-simple completely regular Banach-algebras
 $A$ (of which $C(X)$ is an example) with the
integers under an automorphism of $A$. In particular, various
 properties
of the ideal structure in such crossed products are shown to be
 equivalent to topological properties
of the naturally induced topological dynamical system on $\Delta(A)$,
 the character space of $A$.

In this paper, we recall some of the most important results from
 \cite{svesildej} and \cite{art2}, and in a number of cases provide
significantly simplified proofs of generalizations of results
 occurring
 in \cite{art2}, giving a clearer view of the heart of the matter.
We also include results of a new type in the algebraic setup, and
 furthermore start the investigation of the Banach algebra crossed
product $\ell_1^\sigma (\mathbb{Z}, A)$ of a commutative $C^*$-algebra
 $A$ with the integers under an automorphism $\sigma$ of $A$. In the
 case when $A$ is unital,
this algebra is precisely the one whose $C^*$-envelope is the crossed
 product $C^*$-algebra mentioned above.

This paper is organized as follows.
In Section~\ref{defbas} we give the most general definition of the
kind of crossed product that we will use throughout the first sections
 of this paper. We
also mention the elementary result that the commutant of the
coefficient algebra is automatically a maximal commutative subalgebra
 of the
crossed product.

In Section~\ref{comint} we prove that for any such crossed product $A
 \rtimes_{\Psi} \mathbb{Z}$,
the commutant $A'$ of the coefficient algebra $A$ has non-zero
 intersection with any non-zero ideal $I \subseteq A \rtimes_{\Psi}
 \mathbb{Z}$.
In \cite[Theorem 6.1]{art2}, a more complicated proof of this was
 given
 for a restricted class of coefficient algebras $A$.

In Section~\ref{auto} we focus on the case when $A$ is a function
 algebra on a set $X$ with an automorphism $\widetilde{\sigma}$ of $A$
 induced by a bijection $\sigma : X \rightarrow X$.
According to \cite[Theorem 5.4]{tom2}, the following three
properties are equivalent for a compact Hausdorff space $X$ and a
 homeomorphism $\sigma$ of $X$:
\begin{itemize}
\item The non-periodic points of $(X, \sigma)$ are dense in $X$;
\item Every non-zero closed ideal $I$ of the crossed product
 $C^*$-algebra $C(X) \rtimes_{\alpha} \mathbb{Z}$ is such that $I \cap
 C(X) \neq
 \{0\}$;
\item $C(X)$ is a maximal abelian $C^*$-subalgebra of $C(X)
 \rtimes_{\alpha} \mathbb{Z}$.
\end{itemize}

In Theorem~\ref{bairemax} an analogue of this result is proved for
 our
setup. A reader familiar with the theory of crossed product
$C^*$-algebras will easily recognize that if one chooses $A = C(X)$
for $X$ a compact Hausdorff space in this theorem,
then the algebraic crossed product is canonically isomorphic to a
 norm-dense
subalgebra of the crossed product $C^*$-algebra coming from the
considered induced dynamical system.

For a different kind coefficient algebras $A$ than the ones allowed in
 Theorem~\ref{bairemax},
we prove a similar result in Theorem~\ref{domunmax}.
Theorem~\ref{bairemax} and Theorem~\ref{domunmax} have no non-trivial
 situations in common (Remark~\ref{disrmk}).

In Section~\ref{between} we show that in many situations we can
 always find both a subalgebra properly between the coefficient
 algebra
 $A$
 and its commutant $A'$ (as long as $A \subsetneq A'$, a property we
 have
 a precise condition for in Theorem~\ref{bairemax}) and a non-trivial
 ideal trivially intersecting it, and a subalgebra properly between $A$
 and $A'$ intersecting every non-trival ideal non-trivially.

Section~\ref{setbas} is concerned with the algebraic crossed product
 of
 a complex commutative semi-simple Banach algebra $A$ with the
 integers
 under an automorphism $\sigma$ of $A$, naturally inducing a
homeomorphism $\widetilde{\sigma}$ of the character space $\Delta(A)$
 of $A$. We extend results from~\cite{art2}.

In Section~\ref{bac} we introduce the Banach algebra crossed product
 $\ell_1^{\sigma} (\mathbb{Z}, A)$ for a commutative $C^*$-algebra $A$
 and
 an automorphism $\sigma$ of $A$.
In Theorem~\ref{commker} we give an explicit description of the closed
 commutator ideal in this algebra in terms of the dynamical system
 naturally induced on $\Delta(A)$. We determine the characters of
 $\ell_1^{\sigma} (\mathbb{Z}, A)$. The modular ideals which are
 maximal and
 contain the commutator ideal are precisely the kernels of the
 characters.
\section{Definition and a basic result}\label{defbas}
Let $A$ be an associative commutative complex algebra and let $\Psi
: A \rightarrow A$ be an algebra automorphism. Consider the set
\[A \rtimes_{\Psi} \mathbb{Z} = \{f: \mathbb{Z} \rightarrow A \,|\,
 f(n) = 0 \,\,\textup{except for a finite number of}\, \,n\}.\]
We endow it with the structure of an associative complex algebra by
defining scalar multiplication and addition as the usual pointwise
operations. Multiplication is defined by \emph{twisted convolution},
$*$, as follows;
\[(f*g) (n) = \sum_{k \in \mathbb{Z}} f(k) \cdot \Psi^k (g(n-k)),\]
where $\Psi^k$ denotes the $k$-fold composition of $\Psi$ with
itself. It is trivially verified that $A \rtimes_{\Psi}
\mathbb{Z}$ \emph{is} an associative $\mathbb{C}$-algebra under
these operations. We call it the \emph{crossed product of $A$ and
$\mathbb{Z}$ under $\Psi$}.

A useful way of working with $A \rtimes_{\Psi} \mathbb{Z}$ is to
write elements $f, g \in A \rtimes_{\Psi} \mathbb{Z}$ in the form $f
= \sum_{n \in \mathbb{Z}} f_n \delta^n, g = \sum_{m \in \mathbb{Z}}
g_m \delta^n$, where $f_n = f(n), g_m = g(m)$, addition and scalar
multiplication are canonically defined, and multiplication is
determined by $(f_n \delta^n)*(g_m \delta^m) = f_n \cdot \Psi^n
(g_m) \delta^{n+m}$, where $n,m \in \mathbb{Z}$ and $f_n, g_m \in A$
are arbitrary.

Clearly one may canonically view $A$ as an abelian subalgebra of $A
\rtimes_{\Psi} \mathbb{Z}$, namely as $\{f_0 \delta^0 \, | \, f_0
\in A\}$. The following elementary result is proved in
\cite[Proposition 2.1]{svesildej}.

\begin{pro}\label{comcom}
The commutant $A'$ of $A$ is abelian, and thus it is the unique
maximal abelian subalgebra containing $A$.
\end{pro}
\section{Every non-zero ideal has non-zero intersection with
 $A'$}\label{comint}
Throughout the whole paper, when speaking of an ideal we shall always
 mean a \emph{two-sided} ideal.
We shall now show that any non-zero ideal in $A \rtimes_{\Psi}
\mathbb{Z}$ has non-zero intersection with $A'$. This result,
 Theorem~\ref{commint}, should be
compared with Theorem~\ref{bairemax}, which says that a non-zero ideal
 may well intersect $A$ solely in $0$. An analogue of Theorem~\ref{commint} in the context of
 crossed product $C^*$-algebras is found in \cite[Theorem 4.3]{junojag}. 
Note that in~\cite{art2} a proof of Theorem~\ref{commint} was given
 for
 the case when $A$ was completely regular semi-simple Banach algebra,
 and that this
proof heavily relied upon $A$ having these properties. The present
 proof is elementary and valid for arbitrary commutative algebras.
\begin{thm}\label{commint}
Let $A$ be an associative commutative complex algebra, and let\\
$\Psi : A \rightarrow A$ be an automorphism. Then every non-zero ideal
 of $A \rtimes_{\Psi} \mathbb{Z}$ has non-zero
intersection with the commutant $A'$ of $A$.
\end{thm}
\begin{proof}
Let $I$ be a non-zero ideal, and let $f = \sum_n f_n \delta^n \in I$
 be
 non-zero.
Suppose that $f \notin A'$. Then there must be an $f_{n_i}$ and $a \in
 A$ such that $f_{n_i} \cdot a \neq 0$.
Hence $f' := (\sum_{n} f_n \delta^n) * \Psi^{-n_i} (a) \delta^{-n_i}$
 is a non-zero element of $I$, having $f_{n_i} \cdot a$ as coefficient
 of
 $\delta^0$ and having at most as many non-zero coefficients as $f$.
If $f' \in A'$ we are done, so assume $f' \notin A'$.
Then there exists $b \in A$ such that $F:= b*f' - f'*b \neq 0$.
 Clearly
 $F \in I$ and it is easy to see that $F$ has strictly less non-zero
 coefficients than $f'$ (the coefficient of $\delta^0$ in $F$ is
 zero),
 hence strictly less than $f$. Now if $F \in A'$, we are done. If not,
 we
 repeat the above procedure. Ultimately, if we do not happen to obtain
 a
 non-zero element of $I \cap A'$ along the way, we will be left with a
 non-zero monomial $G := g_m \delta^m \in I$. If this does not lie in
 $A'$, there is an $a\in A$ such that $g_m \cdot a \neq 0$. Hence $G *
 \Psi^{-m} (a) \delta^{-m} = g_m \cdot a \in I \cap A \subseteq I \cap
 A'$.
\end{proof}
Note that the fact that all elements in $A \rtimes_{\Psi} \mathbb{Z}$
 are \emph{finite} sums of the form $\sum_n f_n \delta^n$ is crucial
 for
 the argument used in the proof.
\section{Automorphisms induced by bijections}\label{auto}
Fix a non-empty set $X$, a bijection $\sigma : X \rightarrow X$, and
 an
 algebra of functions $A \subseteq \mathbb{C}^X$ that is invariant
 under $\sigma$ and $\sigma^{-1}$, i.e., such that if $h \in A$, then
 $h
 \circ \sigma \in A$ and $h \circ \sigma^{-1} \in A$.
Then $(X, \sigma)$ is a discrete dynamical system (the action of $n
 \in
 \mathbb{Z}$ on $x \in X$ is given by $n: x \mapsto \sigma^n (x)$) and
 $\sigma$ induces an automorphism $\widetilde{\sigma} : A \rightarrow
 A$
 defined by $\widetilde{\sigma} (f) = f \circ \sigma^{-1}$ by which
 $\mathbb{Z}$ acts on $A$ via iterations.

In this section we will consider the crossed product $A
 \rtimes_{\widetilde{\sigma}} \mathbb{Z}$ for the above setup, and
 explicitly describe
 the commutant $A'$ of $A$. Furthermore, we will investigate
 equivalences between properties of non-periodic points of the system
 $(X,
 \sigma)$, and properties of $A'$.
First we make a few definitions.
\begin{defn}\label{sepper}
For any nonzero $n \in \mathbb{Z}$ we set
\begin{align*}
\sep_A^n (X) &= \{x \in X | \exists h \in A :  h(x) \neq
 \widetilde{\sigma}^n (h) (x)\}, \\
\Per_A^n (X) &= \{x \in X | \forall h \in A : h(x) =
 \widetilde{\sigma}^n (h) (x)\},\\
\sep^n (X) &= \{x \in X |x \neq \sigma^{n} (x))\}, \\
\Per^n (X) &= \{x \in X |x = \sigma^{n} (x)\}.
 \end{align*}
Furthermore, let
\begin{align*}
\Per_A^{\infty} (X) &= \bigcap_{n \in \mathbb{Z} \setminus \{0\}}
 \sep_A^n (X),\\
\Per^{\infty}(X) &= \bigcap_{n \in \mathbb{Z} \setminus \{0\}} \sep^n
 (X).\\
\end{align*}
Finally, for $f \in A$, put
\begin{align*}\supp(f) &= \{x \in X \, | \, f(x) \neq 0\}.
\end{align*}
\end{defn}
It is easy to check that all these sets, except for $\supp(f)$, are
 $\mathbb{Z}$-invariant and that if $A$ separates the points of $X$,
 then
 $\sep_A^n (X) = \sep^n (X)$ and $\Per_A^n (X) = \Per^n(X)$.
Note also that $X \setminus \Per_A^n (X) = \sep_A^n (X)$, and $X
 \setminus \Per^n (X) = \sep^n (X)$. Furthermore $\sep_A^n (X) =
 \sep_A^{-n}
 (X)$ with similar equalities for $n$ and $-n$ ($n \in \mathbb{Z}$)
 holding for $\Per_A^n (X)$, $\sep^n (X)$ and $\Per^n (X)$ as well.
\begin{defn}\label{domun}
We say that a non-empty subset of $X$ is a \emph{domain of uniqueness
 for $A$} if every function in $A$ that vanishes on it, vanishes on
 the
 whole of $X$.
\end{defn}
For example, using results from elementary topology one easily shows
 that for a completely regular topological space $X$, a subset of $X$
 is a
 domain of uniqueness for $C(X)$ if and only if it is dense in $X$.
In the following theorem we recall some elementary results from
 \cite{svesildej}.
\begin{thm}\label{commdesc}
The unique maximal abelian subalgebra of $A
 \rtimes_{\widetilde{\sigma}} \mathbb{Z}$ that contains $A$ is
 precisely the set of elements
\[A' =\{\sum_{n \in \mathbb{Z}} f_n \delta^n \,|\, f_n
 \!\!\upharpoonright_{\sep^n_A (X)} \equiv 0 \textup{ for all } n \in
 \mathbb{Z}\}.\]
So if $A$ separates the points of $X$, then
\[A' =\{\sum_{n \in \mathbb{Z}} f_n \delta^n \,|\, \supp(f_n)
 \subseteq
 \Per^n (X) \textup{ for all } n \in \mathbb{Z}\}.\]
Furthermore, the subalgebra $A$ is maximal abelian in $A
 \rtimes_{\widetilde{\sigma}} \mathbb{Z}$ if and only if, for every $n
 \in \mathbb{Z}
 \setminus \{0\}$, $\sep_A^n (X)$ is a domain of uniqueness for $A$.
\end{thm}
We now focus solely on topological contexts. In order to prove one of
 the main theorems of this section, we need the following topological
 lemma.
\begin{lem}\label{bairegrej}
Let $X$ be a Baire space which is also Hausdorff, and let $\sigma : X
 \rightarrow X$ be a homeomorphism.
Then the non-periodic points of $(X, \sigma)$ are
dense if and only if $\Per^{n} (X)$ has empty interior for all
 positive
 integers $n$.
\end{lem}
\begin{proof}
Clearly, if there is a positive integer $n$ such that $\Per^{n}
(X)$ has non-empty interior, the non-periodic points are not
dense. For the converse we note that we may write
\[X \backslash \Per^{\infty} (X) = \bigcup_{n > 0} \Per^n (X).\]
If the set of non-periodic points is not dense, its complement has
non-empty interior, and as the sets $\Per^n (\Delta(A))$ are clearly
all closed since $X$ is Hausdorff, there must exist an integer $n_0 >
 0$ such that
$\Per^{n_0} (X)$ has non-empty interior since $X$ is
a Baire space. \end{proof}
We are now ready to prove the following theorem.
\begin{thm}\label{bairemax}
Let $X$ be a Baire space which is also Hausdorff, and let $\sigma : X
 \rightarrow X$ be a homeomorphism inducing, as usual, an automorphism
 $\widetilde{\sigma}$ of $C(X)$.  Suppose $A$ is a subalgebra of
 $C(X)$
 that is invariant under $\widetilde{\sigma}$ and its inverse,
 separates
 the points of $X$ and is such that for every non-empty open set $U
 \subseteq X$ there is a non-zero $f \in A$ that vanishes on the
 complement
 of $U$. Then the following three statements are equivalent.

\begin{itemize}
\item $A$ is a maximal abelian subalgebra of  $A
 \rtimes_{\widetilde{\sigma}} \mathbb{Z}$;
\item $\Per^{\infty} (X)$ is dense in $X$;
\item Every non-zero ideal $I \subseteq A \rtimes_{\widetilde{\sigma}}
 \mathbb{Z}$ is such that $I \cap A \neq \{0\}$.
\end{itemize}
\end{thm}
\begin{proof}
Equivalence of the first two statements is precisely the result in
 \cite[Theorem 3.7]{svesildej}. The first property implies the third
 by
 Proposition~\ref{comcom} and Theorem~\ref{commint}. Finally, to show
 that
 the third statement implies the second, assume that $\Per^\infty (X)$
 is
 not dense.
It follows from Lemma~\ref{bairegrej} that there exists an integer $n
 >0$ such
that $\Per^n (X)$ has non-empty interior. By the assumptions on $A$
 there exists a non-zero $f \in A$ such that $\supp(f)
\subseteq \Per^n (X)$. Consider now the non-zero ideal $I$ generated
 by
 $f + f
\delta^n$. It is spanned by elements of the form $a_i
\delta^i * (f + f \delta^n) * a_j \delta^j$, $(f+f\delta^n) * a_j
 \delta^j$, $a_i \delta^i * (f + f \delta^n)$ and $f+f\delta^n$. Using
 that
 $f$ vanishes
outside $\Per^n (X)$, so that $f \delta^n * a_j \delta^j = a_j f
 \delta^{n+j}$, we may for example rewrite
\begin{align*}
&a_i \delta^i * (f + f \delta^n)* a_j \delta^j = [a_i \cdot (a_j \circ
\widetilde{\sigma}^{-i})\delta^{i}]* [f \delta^j + f \delta^{n+j}]\\
&= [a_i \cdot (a_j \circ \widetilde{\sigma}^{-i}) \cdot (f\circ
\widetilde{\sigma}^{-i})] \delta^{i+j} + [a_i \cdot (a_j \circ
\widetilde{\sigma}^{-i}) \cdot (f\circ \widetilde{\sigma}^{-i})]
\delta^{i+j+n}.
\end{align*}
A similar calculation for the other three kinds of elements that span
 $I$ now makes it clear that any element in $I$ may be written in the
 form
$\sum_{i} (b_i \delta^i + b_i\delta^{n+i})$. As $i$ runs only
through a finite subset of $\mathbb{Z}$, this is not a non-zero
monomial. In particular, it is not a non-zero element in
$A$. Hence $I$ intersects $A$ trivially.
\end{proof}
We also have the following result for a different kind of subalgebras
 of $C(X)$.
\begin{thm}\label{domunmax}
Let $X$ be a topological space, $\sigma : X \rightarrow X$ a
 homeomorphism, and $A$ a non-zero subalgebra of $C(X)$, invariant
 both
 under the
 usual induced automorphism $\widetilde{\sigma}: C(X) \rightarrow
 C(X)$
 and under its inverse. Assume that $A$ separates the points of $X$
 and
 is such that every non-empty open set $U \subseteq X$ is a domain of
 uniqueness for $A$. Then the following three statements are
 equivalent.
\begin{itemize}
\item $A$ is maximal abelian in $A \rtimes_{\widetilde{\sigma}}
 \mathbb{Z}$;
\item $\sigma$ is not of finite order;
\item Every non-zero ideal $I \subseteq A \rtimes_{\widetilde{\sigma}}
 \mathbb{Z}$ is such that $I \cap A \neq \{0\}$.
\end{itemize}
\end{thm}
\begin{proof}
Equivalence of the first two statements is precisely the result in
 \cite[Theorem 3.11]{svesildej}.
That the first statement implies the third follows immediately from
 Proposition~\ref{comcom} and Theorem~\ref{commint}.
Finally, to show that the third statement implies the second, assume
 that there exists an $n$, which we may clearly choose to be
 non-negative,
 such that $\sigma^n = \id_X$. Now take any non-zero $f \in A$ and
 consider the non-zero ideal $I = (f + f \delta^n)$. Using an argument
 similar to the one in the proof of Theorem~\ref{bairemax} one
 concludes
 that $I \cap A = \{0\}$.
\end{proof}
\begin{cor}\label{comfd}
Let $M$ be a connected complex manifold and suppose the function
 $\sigma : M \rightarrow M$ is biholomorphic. If $A \subseteq H(M)$ is
 a
 subalgebra of the algebra of holomorphic functions that separates the
 points
 of $M$ and which is invariant under the induced automorphism
 $\widetilde{\sigma}$ of $H(M)$ and its inverse, then the following
 three
 statements are equivalent.
\begin{itemize}
\item $A$ is maximal abelian in $A \rtimes_{\widetilde{\sigma}}
 \mathbb{Z}$;
\item $\sigma$ is not of finite order;
\item Every non-zero ideal $I \subseteq A \rtimes_{\widetilde{\sigma}}
 \mathbb{Z}$ is such that $I \cap A \neq \{0\}$.
\end{itemize}
\end{cor}
\begin{proof}
On connected complex manifolds, open sets are domains of uniqueness
 for
 $H(M)$. See for example \cite{complex}.
\end{proof}
\begin{rmk}\label{disrmk}
It is worth mentioning that the required conditions in
 Theorem~\ref{bairemax} and Theorem~\ref{domunmax} can only be
 simultaneously
 satisfied in case $X$ consists of a single point and $A =
 \mathbb{C}$.
 This
 is explained in \cite[Remark 3.13]{svesildej}
\end{rmk}
\section{Algebras properly between the coefficient algebra and its
 commutant}\label{between}
From Theorem~\ref{bairemax} it is clear that for spaces $X$ which are
 Baire and Hausdorff and subalgebras $A \subseteq C(X)$ with
 sufficient
 separation properties, $A$ is equal to its own commutant in the
 associated crossed product precisely when the aperiodic points,
 $\Per^\infty
 (X)$, constitute a dense subset of $X$. This theorem also tells us
 that
 whenever $\Per^\infty (X)$ is not dense there exists a non-zero
 ideal
 $I$ having zero intersection with $A$, while the general
 Theorem~\ref{commint} tells us that every non-zero ideal has non-zero
 intersection
 with $A'$, regardless of the system $(X, \sigma)$.
\begin{defn}
We say that a subalgebra has the \emph{intersection property} if it
 intersects every non-zero ideal non-trivially.
\end{defn}
A subalgebra $B$ such that $A \subsetneq B \subsetneq A'$ is said to be properly between $A$ and $A'$.
Two natural questions comes to mind in case $\Per^{\infty}(X)$ is not
 dense:
\begin{enumerate}
 \item Do there exist subalgebras properly between $A$ and $A'$ having
 the intersection property?
\item Do there exist subalgebras properly between $A$ and $A'$
 \emph{not} having the intersection property?
\end{enumerate}
 We shall show that for a significant class of systems the answer to
 both these questions is positive.
\begin{pro}\label{intalgid}
Let $X$ be a Hausdorff space, and let $\sigma : X
 \rightarrow X$ be a homeomorphism inducing, as usual, an automorphism
 $\widetilde{\sigma}$ of $C(X)$.  Suppose $A$ is a subalgebra of
 $C(X)$
 that is invariant under $\widetilde{\sigma}$ and its inverse,
 separates
 the points of $X$ and is such that for every non-empty open set $U
 \subseteq X$ there is a non-zero $f \in A$ that vanishes on the
 complement
 of $U$. Suppose furthermore that there exists an integer $n > 0$ such
 that the interior of $\Per^n (X)$ contains at least two orbits. Then
 there exists a subalgebra $B$ such that $A \subsetneq B \subsetneq
 A'$ which does not have the intersection property.
\end{pro}
\begin{proof}
Using the Hausdorff property of $X$ and the fact that $\Per^n (X)$
 contains two orbits we can find two non-empty disjoint invariant open
 subsets $U_1$ and $U_2$ contained in $\Per^n (X)$.
Consider
\[B = \{f_0 + \sum_{k \neq 0} f_k \delta^k : f_0 \in A, \, \supp(f_k)
 \subseteq U_1 \cap \Per^k (X) \textup{ for } k \neq 0\}.\]
Then $B$ is a subalgebra and $B \subseteq A'$.
The assumptions on $A$ and the definitions of $U_1$ and $U_2$ now
 make it clear that $A \subsetneq B \subsetneq A'$ since there exist,
 for
 example, non-zero functions $F_1, F_2 \in A$ such that $\supp (F_1)
 \subseteq U_1$ and $\supp (F_2) \subseteq U_2$, and thus $F_1
 \delta^n
 \in B
 \setminus A$ and $F_2 \delta^n \in A' \setminus B$.
Consider the non-zero ideal $I$ generated by $F_2 + F_2 \delta^n$.
 Using an
 argument similar to the one used in the proof of
 Theorem~\ref{bairemax} we
 see that $I \cap A = \{0\}$. It is also easy to see that $I \subseteq
 \{\sum_k f_k \delta^k : \supp(f_k) \subseteq U_2\}$ since $U_2$ is
 invariant. As $U_1 \cap U_2 = \emptyset$, we see from the description
 of $B$
 that $I \cap B \subseteq A$, so that $I \cap B \subseteq I \cap A =
 \{0\}$.
\end{proof}
We now exhibit algebras properly between $A$ and $A'$ that do have the intersection property.
\begin{pro}\label{intprop}
Let $X$ be a Hausdorff space, and let $\sigma : X
 \rightarrow X$ be a homeomorphism inducing, as usual, an automorphism
 $\widetilde{\sigma}$ of $C(X)$.  Suppose $A$ is a subalgebra of
 $C(X)$
 that is invariant under $\widetilde{\sigma}$ and its inverse,
 separates
 the points of $X$ and is such that for every non-empty open set $U
 \subseteq X$ there is a non-zero $f \in A$ that vanishes on the
 complement
 of $U$. Suppose furthermore that there exist an integer $n > 0$ such
 that the interior of $\Per^n (X)$ contains a point $x_0$ which is not
 isolated, and an $f \in A$ with $\supp(f) \subseteq \Per^n (X)$ and
 $f(x_0) \neq 0$. Then
 there exists a subalgebra $B$ such that $A \subsetneq B \subsetneq
 A'$
 which has the intersection property.
\end{pro}
\begin{proof}
Define
\[B = \{\sum_{k \in \mathbb{Z}} f_k \delta^k \in A' \,|\, f_k (x_0) =
 0
 \textup{ for all } k \neq 0\},\]
where $x_0$ is as in the statement of the theorem. Clearly $B$ is a
 subalgebra
and $A \subseteq B$. Since $x_0$ is not isolated, we can use the
 assumptions on
$A$ and the fact that $X$ is Hausdorff to first find a point different
 from
$x_0$ in the interior of $\Per^n (X)$ and subsequently a non-zero
 function $g
\in A$ such that $\supp(g) \subseteq \Per^n (X)$ and $g(x_0) = 0$. Then
 $g
\delta^n \in B \setminus A$. Also, by the assumptions on $A$ there is
 a
non-zero $f \in A$ with $\supp(f) \subseteq \Per^n (X)$ such that
 $f(x_0) \neq
0$, whence $f \delta^n \in A' \setminus B$. This shows that $B$ is a
 subalgebra
properly between $A$ and
 $A'$.
To see that it has the intersection property, let $I$ be an arbitrary
 non-zero ideal in the crossed product and note that by
 Theorem~\ref{commint} there is a non-zero $F = \sum_{k \in \mathbb{Z}}
 f_k \delta^k$ in $I \cap A'$.
Now if for all $k \neq 0$ we have that $f_k (x_0) = 0$, we are done.
 So
 suppose there is some $k \neq 0$ such that $f_k (x_0) \neq 0$. Since
 $f_k$ is continuous and $x_0$ is not isolated, we may use the
 Hausdorff property of $X$ to conclude that there exists a non-empty open set
 $V$ contained in the interior of $\Per^n (X)$ such that $x_0 \notin V$
 and $f_k (x) \neq 0$ for all $x \in V$. The
 assumptions on $A$ now imply that there is an $h \in A$ such that
 $h(x_0) = 0$ and $h(x_1) \neq 0$ for some $x_1 \in V \subseteq
 \supp(f_n)$. Clearly $0 \neq h * F \in I \cap B$.
\end{proof}
\begin{thm}\label{corsnitt}
Let $X$ be a Baire space which is Hausdorff and connected. Let $\sigma:
 X \rightarrow X$ be a
 homeomorphism inducing an automorphism $\widetilde{\sigma}$ of $C(X)$
 in the
 usual way.  Suppose $A$ is a subalgebra of $C(X)$ that is invariant
 under $\widetilde{\sigma}$ and its inverse, such that for every open
 set $U \subseteq X$ and $x \in U$ there is an $f \in A$ such that $f (x)
 \neq 0$ and $\supp(f) \subseteq U$.
Then precisely one of the following situations occurs:
\begin{enumerate}
\item $A = A'$, which happens precisely when $\Per^{\infty} (X)$ is
 dense;
\item $A \subsetneq A'$ and there exist both subalgebras properly
 between $A$
and $A'$ which have the intersection property, and subalgebras which do
 not.
This happens precisely when $\Per^{\infty} (X)$ is not dense and $X$
 is
infinite; \item $A \subsetneq A'$ and every subalgebra properly between
 $A$ and
$A'$ has the intersection property. This happens precisely when $X$
 consists of
one point.
\end{enumerate}
\end{thm}
\begin{proof}
By Theorem~\ref{bairemax}, (i) is clear and we may assume that
 $\Per^{\infty} (X)$ is not dense. Suppose first that $X$ is infinite and note
 that by
 Lemma~\ref{bairegrej} there exists $n_0 > 0$ such that $\Per^{n_0}
 (X)$ has
 non-empty interior. If this interior consists of one single orbit
 then as $X$ is Hausdorff every point in the interior is both closed
 and
 open, so that $X$ consists of one point by connectedness, which is a
 contradiction. Hence there are at least two orbits in the interior of
 $\Per^{n_0} (X)$. Furthermore, no point of $X$ can be isolated. Thus
 by Proposition~\ref{intalgid} and Proposition~\ref{intprop} there are
 subalgebras properly between $A$ and $A'$ which have the intersection
 property, and subalgebras which do not.
Suppose next that $X$ is finite, so that $X = \{x\}$ by connectedness.
 Then
$\sigma$ is the identity map, and $A = \mathbb{C}$. In this case, $A
\rtimes_{\widetilde{\sigma}} \mathbb{Z}$ may be canonically identified
 with
$\mathbb{C} [t, t^{-1}]$. Let $B$ be a subalgebra such that
 $\mathbb{C}
\subsetneq B \subsetneq \mathbb{C}[t, t^{-1}]$, and let $I$ be a
 non-zero ideal
of $\mathbb{C}[t, t^{-1}]$. We will show that $I\cap B\neq \{0\}$ and
 hence may
assume that $I\neq \mathbb{C}[t, t^{-1}]$. Since $\mathbb{C} [t,
 t^{-1}]$ is
the ring of fractions of $\mathbb{C}[t]$ with respect to the
 multiplicatively
closed subset $\{t^n \,|\, n \textup{ is a non-negative integer}\}$
 and
$\mathbb{C}[t]$ is a principal ideal domain, it follows
 from~\cite[Proposition
3.11, (i)]{atmac} that $I$ is of the form $(t-\alpha_1)\cdots
 (t-\alpha_n)
\mathbb{C}[t,t^{-1}]$ for some $n>0$ and $\alpha_1, \ldots, \alpha_n
 \in
\mathbb{C}$. There exists a non-constant $f$ in $B$, and then the
 element
$(f-f(\alpha_1)) \cdots (f - f(\alpha_n))$ is a non-zero element of
 $B$. It is
clearly also in $I$ since it vanishes at $\alpha_1, \ldots, \alpha_n$ and hence has $(t-\alpha_1)\cdots(t-\alpha_n)$ as a factor. Hence $I \cap B \neq \{0\}$ and the proof is completed.
\end{proof}
It is interesting to mention that arguments similar to the ones used in
 Propositions~\ref{intalgid} and~\ref{intprop} work in the
 context of the crossed product $C^*$-algebra $C(X) \rtimes_{\alpha}
 \mathbb{Z}$ where $X$ is a compact Hausdorff space and $\alpha$ the
 automorphism induced by a homeomorphism of $X$. See \cite[Section 5]{junojag} for details.
 \section{Semi-simple Banach algebras}\label{setbas}
In what follows, we shall focus on cases where $A$ is a commutative
complex Banach algebra, and freely make use of the basic theory for
such $A$, see e.g.\!\! \cite{larsen}. As conventions tend to differ
slightly in the literature, however, we mention that we call a
commutative Banach algebra $A$ \emph{completely regular} (the term
 \emph{regular} is also frequently used in the literature)
if, for every subset $F \subseteq \Delta(A)$ (where $\Delta(A)$
 denotes
 the character space of $A$) that is closed in the
Gelfand topology and for every $\phi_0 \in \Delta(A) \setminus F$,
there exists an $a \in A$ such that $\widehat{a} (\phi) = 0$ for all
$\phi \in F$ and $\widehat{a} (\phi_0) \neq 0$. All topological
considerations of $\Delta(A)$ will be done with respect to its
Gelfand topology.

Now let $A$ be a complex commutative semi-simple completely regular
 Banach
algebra, and let $\sigma  : A \rightarrow A$ be an algebra
automorphism. As in \cite{svesildej}, $\sigma$ induces a map
$\widetilde{\sigma}: \Delta(A) \rightarrow \Delta(A)$ defined by
$\widetilde{\sigma}(\mu) = \mu \circ \sigma^{-1}, \mu \in
\Delta(A)$, which is automatically a homeomorphism when $\Delta(A)$
is endowed with the Gelfand topology. Hence we obtain a topological
dynamical system $(\Delta(A), \widetilde{\sigma})$. In turn,
$\widetilde{\sigma}$ induces an automorphism $\widehat{\sigma} :
\widehat{A} \rightarrow \widehat{A}$ (where $\widehat{A}$ denotes
the algebra of Gelfand transforms of all elements of $A$) defined by
$\widehat{\sigma} (\widehat{a}) = \widehat{a} \circ
\widetilde{\sigma}^{-1} = \widehat{\sigma(a)}$. Therefore we can
form the crossed product $\widehat{A} \rtimes_{\widehat{\sigma}}
\mathbb{Z}$.

In what follows, we shall make frequent use of the following fact.
Its proof consists of a trivial direct verification.
\begin{thm}\label{isom}
Let $A$ be a commutative semi-simple Banach algebra and $\sigma$ an
automorphism, inducing an automorphism $\widehat{\sigma} :
\widehat{A} \rightarrow \widehat{A}$ as above. Then the map $\Phi :
A \rtimes_{\sigma} \mathbb{Z} \rightarrow \widehat{A}
\rtimes_{\widehat{\sigma}} \mathbb{Z}$ defined by $\sum_{n \in
\mathbb{Z}} a_n \delta^n \mapsto \sum_{n \in \mathbb{Z}}
\widehat{a_n} \delta^n$ is an isomorphism of algebras mapping $A$
onto $\widehat{A}$.
\end{thm}
We shall now conclude that, for certain $A$, two
different algebraic properties of $A \rtimes_{\sigma} \mathbb{Z}$
are equivalent to density of the non-periodic points of the
naturally associated dynamical system on the character space
$\Delta(A)$. The analogue of this result in the context of crossed
product $C^*$-algebras is \cite[Theorem 5.4]{tom2}. We shall also
combine this with a theorem from \cite{svesildej} to conclude a
stronger result for the Banach algebra $L_1 (G)$, where $G$ is a
locally compact abelian group with connected dual group.

\begin{thm}\label{triquiv}
Let $A$ be a complex commutative semi-simple completely regular Banach
 algebra,
$\sigma: A \rightarrow A$ an automorphism and $\widetilde{\sigma}$
the homeomorphism of $\Delta(A)$ in the Gelfand topology induced by
$\sigma$ as described above. Then the following three properties are
equivalent:
\begin{itemize}
\item The non-periodic points $\Per^\infty (\Delta(A))$ of
 $(\Delta(A),\widetilde{\sigma})$ are dense in $\Delta(A)$;
\item Every non-zero ideal $I \subseteq A \rtimes_{\sigma} \mathbb{Z}$
 is such that $I \cap A \neq \{0\}$;
\item $A$ is a maximal abelian subalgebra of $A \rtimes_{\sigma}
 \mathbb{Z}$.
\end{itemize}
\end{thm}
\begin{proof}
As $A$ is completely regular, and $\Delta(A)$ is Baire since it is
 locally compact and Hausdorff, it is immediate from
 Theorem~\ref{bairemax}
 that the following three statements are equivalent.
\begin{itemize}
\item The non-periodic points $\Per^\infty (\Delta(A))$ of
 $(\Delta(A),\widetilde{\sigma})$ are dense in $\Delta(A)$;
\item Every non-zero ideal $I \subseteq \widehat{A}
 \rtimes_{\widehat{\sigma}} \mathbb{Z}$ is such that $I \cap
 \widehat{A} \neq \{0\}$;
\item $\widehat{A}$ is a maximal abelian subalgebra of $\widehat{A}
 \rtimes_{\widehat{\sigma}} \mathbb{Z}$.
\end{itemize}
Now applying Theorem~\ref{isom} we can pull everything back to $A
 \rtimes_{\sigma} \mathbb{Z}$ and the result follows.
\end{proof}

The following result for a more specific class of Banach algebras is
 an
 immediate consequence of Theorem~\ref{triquiv} together with
 \cite[Theorem 4.16]{svesildej}.
\begin{thm}\label{conndualtri}
Let $G$ be a locally compact abelian group with connected dual group
and let $\sigma : L_1 (G) \rightarrow L_1 (G)$ be an automorphism.
Then the following three statements are equivalent.
\begin{itemize}
\item $\sigma$ is not of finite order;
\item Every non-zero ideal $I \subseteq L_1 (G) \rtimes_{\sigma}
 \mathbb{Z}$ is such that $I \cap L_1 (G) \neq \{0\}$;
\item $L_1 (G)$ is a maximal abelian subalgebra of $L_1 (G)
 \rtimes_{\sigma} \mathbb{Z}$.
\end{itemize}
\end{thm}

To give a more complete picture, we also include the results
 \cite[Theorem 5.1]{art2} and \cite[Theorem 7.6]{art2}.
\begin{thm}\label{pritop}
Let $A$ be a complex commutative semi-simple completely regular unital
 Banach
algebra such that $\Delta(A)$ consists of infinitely many points,
and let $\sigma$ be an automorphism of $A$. Then
\begin{itemize}
\item $A \rtimes_{\sigma} \mathbb{Z}$ is simple if and only if the
 associated system
$(\Delta(A), \widetilde{\sigma})$ on the character space is minimal.
\item $A \rtimes_{\sigma} \mathbb{Z}$ is prime if and only if
 $(\Delta(A), \widetilde{\sigma})$ is
topologically transitive.
\end{itemize}
\end{thm}

\section{The Banach algebra crossed product $\ell_1^\sigma
 (\mathbb{Z},
 A)$ for a commutative $C^*$-algebra $A$}\label{bac}
Let $A$ be a commutative $C^*$-algebra with spectrum $\Delta(A)$ and
 $\sigma : A \rightarrow A$ an automorphism.
We identify the set $\ell^1 (\mathbb{Z}, A)$ with the set $\{\sum_{n
 \in \mathbb{Z}} f_n \delta^n | f_n \in A, \sum_{n \in \mathbb{Z}}
 \|f_n\|
 < \infty\}$ and endow it with the same operations as for the finite
 sums in Section~\ref{defbas}.
Using that $\sigma$ is isometric one easily checks that the operations
 are well defined, and that the usual norm on this set is an algebra
 norm with respect to the convolution product.

We denote this algebra by $\ell_1^\sigma (\mathbb{Z}, A)$, and note
 that it is a Banach algebra.
By basic theory of $C^*$-algebras, we have the isometric automorphism
 $A \cong \widehat{A} = C_0 (\Delta(A))$.
As in Section~\ref{setbas}, $\sigma$ induces a homeomorphism,
 $\widetilde{\sigma} : \Delta(A) \rightarrow \Delta(A)$ and an
 automorphism
 $\widehat{\sigma} : C_0(\Delta(A)) \rightarrow C_0 (\Delta(A))$ and
 we
 have
 a canonical isometric isomorphism of $\ell_1^\sigma (\mathbb{Z}, A)$
 onto $\ell_1^{\widehat{\sigma}} (\mathbb{Z}, C_0 (\Delta(A)))$ as in
 Theorem~\ref{isom}.

We will work in the concrete crossed product
 $\ell_1^{\widehat{\sigma}}
 (\mathbb{Z}, C_0 (\Delta(A)))$.
We shall describe the closed commutator ideal $\mathscr{C}$ in terms
 of
 $(\Delta(A), \sigma)$.
In analogy with the notation used in~\cite{hullkertom}, we make the
 following definitions.
\begin{defn}\label{hullker}
Given a subset $S \subseteq \Delta(A)$, we set
\begin{align*}
\ker(S) &= \{f \in C_0 (\Delta(A))\,|\, f(x) = 0 \textup{ for all } x
 \in S\},\\
\Ker(S) &= \{\sum_{n \in \mathbb{Z}} f_n \in
 \ell_1^{\widehat{\sigma}} (\mathbb{Z}, C_0 (\Delta(A)))\,|\, f_n(x) =
 0 \textup{ for all } x
 \in S, n \in \mathbb{Z}\}.
\end{align*}
\end{defn}
Clearly $\Ker(S)$ is always a closed subspace, and in case $S$ is
 invariant, it is a closed ideal.

We shall also need the following version of the Stone-Weierstrass
 theorem.
\begin{thm}\label{stone}
Let $X$ be a locally compact Hausdorff space and let $C$ be a closed
 subset of $X$. Let $B$ be a self-adjoint subalgebra of $C_0 (X)$
 vanishing on $C$. Suppose that for any pair of points $x, y \in X$,
 with $x
 \neq y$, such that at least one of them is not in $C$, there exists
 $f
 \in
 B$ such that $f(x) \neq f(y)$. Then $\overline{B} = \{f \in C_0 (X) :
 f(x) = 0 \textup{ for all } x \in C\}$.
\end{thm}
\begin{proof}
This follows from the more general result \cite[Theorem
 11.1.8]{dixmier}, as it is well
 known that the pure states of $C_0(\Delta(A))$ are precisely the
 point
 evaluations on the locally compact Hausdorff space $\Delta(A)$, and
 that a pure state of a sub-$C^*$-algebra always has a pure state
 extension to the whole $C^*$-algebra. By passing to the one-point
 compactification of $\Delta(A)$, one may also easily derive the result
 from the more
 elementary \cite[Theorem 2.47]{douglas}.
\end{proof}
\begin{defn}\label{appruni}
Let $A$ be a normed algebra. An approximate unit of $A$ is a net
 $\{E_{\lambda}\}_{\lambda \in \Lambda }$ such that for every $a \in
 A$
 we
 have $\lim_{\lambda} \|E_{\lambda} a - a\| = \lim_{\lambda} \|a
 E_{\lambda}  - a\| =0$.
\end{defn}
Recall that any $C^*$-algebra has an approximate unit such that
 $\|E_{\lambda}\| \leq 1$ for all $\lambda \in \Lambda$. In general,
 however,
 an approximate identity need not be bounded.
We are now ready to prove the following result, which is the analogue
 of the first part of \cite[Proposition 4.9]{hullkertom}.
\begin{thm}\label{commker}
$\mathscr{C} = \Ker(\Per^1 (\Delta(A)))$.
\end{thm}
\begin{proof}
It it easily seen that $\mathscr{C} \subseteq \Ker(\Per^1
 (\Delta(A)))$.
For the converse inclusion we choose an approximate identity
 $\{E_{\lambda}\}_{\lambda \in \Lambda}$ for $C_0 (\Delta(A))$ and
 note
 first of
 all that for any $f \in C_0 (\Delta(A))$ we have $f*(E_\lambda
 \delta)
 -
 (E_\lambda \delta)*f = E_\lambda (f - f \circ
 \widetilde{\sigma}^{-1})
 \delta \in \mathscr{C}$. Hence as $\mathscr{C}$ is closed, $(f-f
 \circ
 \widetilde{\sigma}^{-1})\delta \in \mathscr{C}$ for all $f \in C_0
 (\Delta)$. Clearly the set $J = \{g \in C_0 (\Delta(A)) \,|\, g
 \delta
 \in
 \mathscr{C}\}$ is a closed subalgebra (and even an ideal) of $C_0
 (\Delta(A))$. Denote by $I$ the (self-adjoint) ideal of $C_0
 (\Delta(A))$
 generated by the set of elements of the form $f - f \circ
 \widetilde{\sigma}^{-1}$. Note that $I$ vanishes on $\Per^1
 (\Delta(A))$ and that it
 is contained in $J$. Using complete regularity of $C_0 (\Delta(A))$,
 it
 is straightforward to check that for any pair of distinct points $x,
 y
 \in \Delta(A)$, at least one of which is not in $\Per^1 (\Delta(A))$,
 there exists a function $f \in I$ such that $f(x) \neq f(y)$. Hence
 by
 Theorem~\ref{stone} $I$ is dense in $\ker(\Per^1 (\Delta(A)))$, and
 thus
$\{f \delta \,|\, f \in \ker(\Per^1 (\Delta(A)))\} \subseteq
 \mathscr{C}$ since $J$ is closed. So for any $n \in \mathbb{Z}$ and
 $f
 \in
 \ker(\Per^1 (\Delta(A)))$ we have $(f \delta) * (E_\lambda \circ
 \widetilde{\sigma}) \delta^{n-1} = (f E_\lambda) \delta^n \in
 \mathscr{C}$.
 This converges to $f \delta^n$, and hence $\mathscr{C} \supseteq
 \Ker(\Per^1 (\Delta(A)))$.
\end{proof}
Denote the set of non-zero multiplicative linear functionals of
 $\ell_1^{\widehat{\sigma}} (\mathbb{Z}, C_0 (\Delta(A)))$ by $\Xi$.
 We
 shall
 now determine a bijection between $\Xi$ and $\Per^1 (\Delta(A))
 \times
 \mathbb{T}$.
It is a standard result from Banach algebra theory that any $\mu \in
 \Xi$ is bounded and of norm at most one. Since one
 may choose an approximate identity $\{E_{\lambda}\}_{\lambda \in
 \Lambda}$ for $C_0 (\Delta(A))$ such that $\|E_{\lambda}\| \leq 1$
 for
 all
 $\lambda \in \Lambda$ it is also easy to see that $\|\mu\| = 1$.
 Namely,
 given $\mu \in \Xi$ we may choose an $f \in C_0 (\Delta(A))$ such
 that
 $\mu (f) \neq 0$. Then by continuity of $\mu$ we have $\mu (f) =
 \lim_{\lambda} \mu (f E_{\lambda}) = \mu(f) \lim_{\lambda}
 (E_{\lambda})$ and
 hence $\lim_{\lambda} \mu(E_{\lambda}) = 1$.

\begin{lem}\label{apprun}
The limit $\xi := \lim_{\lambda} \mu (E_{\lambda} \delta)$ exists for
 all $\mu \in \Xi$, and is independent of the approximate unit
 $\{E_{\lambda}\}_{\lambda \in \Lambda}$. Furthermore, $\xi \in
 \mathbb{T}$ and
 $\lim_{\lambda} \mu (E_{\lambda} \delta^n) = \xi^n$ for all integers
 $n$.
\end{lem}
\begin{proof}
By continuity and multiplicativity of $\mu$ we have that
 $\lim_{\lambda} \mu(f) \mu  (E_{\lambda} \delta) = \mu(f \delta)$ for
 all $f \in C_0
 (X)$.
So for any $f$ such that $\mu(f) \neq 0$ we have that $\lim_{\lambda}
 \mu(E_{\lambda} \delta) = \frac{\mu(f \delta)}{\mu(f)}$. This shows
 that
 the limit $\xi$ exists and is the same for any approximate unit, and
 using a similar argument one easily sees that $\lim_{\lambda}
 \mu(E_{\lambda} \delta^n)$ also exists and is independent of
 $\{E_{\lambda}\}_{\lambda \in \Lambda}$. For the rest of the proof,
 we
 fix an approximate unit $\{E_{\lambda}\}_{\lambda \in \Lambda}$ such
 that $\|E_{\lambda}\| \leq 1$ for all $\lambda \in \Lambda$. As we
 know
 that $\|\mu\| = 1$, we see that $|\xi|
 \leq 1$.
Now suppose $|\xi| < 1$.
It is easy to see that $\lim_{\lambda} \mu(E_{\lambda}) =1 = \xi^0$.
 Hence also $1 = \lim_{\lambda} \mu (E_{\lambda})^2 = \lim_{\lambda}
 \mu (E_{\lambda}^2)  = \lim_{\lambda}
 \mu((E_{\lambda} \delta)*((E_{\lambda} \circ \widetilde{\sigma})
 \delta^{-1})) = \lim_{\lambda} \mu((E_{\lambda} \delta))\cdot
 \lim_{\lambda} \mu((E_{\lambda}
 \circ \widetilde{\sigma}) \delta^{-1})$.
Now as we assumed $|\xi| <1$, this forces $|\lim_{\lambda} \mu([
 (E_{\lambda} \circ \widetilde{\sigma}) \delta^{-1}])| > 1$, which is
 clearly
 a contradiction since $\|\mu\| = 1$. To prove the last statement we
 note that for any $n$, $\{E_{\lambda} \circ
 \widetilde{\sigma}^{-n}\}_{\lambda \in \Lambda}$ is an approximate
 unit for $C_0 (X)$, and that if
 $\{F_{\lambda}\}_{\lambda \in \Lambda}$ is another approximate unit
 for
 $C_0 (X)$ indexed by the same set $\Lambda$, we have that
 $\{E_{\lambda}
 F_{\lambda}\}_{\lambda \in \Lambda}$ is an approximate unit as well.
 Now note that $\mu(E_{\lambda} \delta)\cdot \mu(E_{\lambda} \delta)
  =
 \mu ((E_{\lambda} \delta)* (E_{\lambda}\delta)) = \mu (E_{\lambda}
 (E_{\lambda} \circ \widetilde{\sigma}^{-1}) \delta^2)$. Using what we
 concluded
 above about independence of approximate units, this shows that $\xi^2
 =
 \lim_{\lambda} \mu (E_{\lambda}\delta)^2 = \lim_{\lambda} \mu
 (E_{\lambda}
 (E_{\lambda} \circ \widetilde{\sigma}^{-1}) \delta^2) =
  \lim_{\lambda} \mu (E_{\lambda}\delta^2)$.
Inductively, we see that $\lim_{\lambda} \mu (E_{\lambda} \delta^n) =
 \xi^n$ for non-negative $n$. As $\mu((E_{\lambda} \delta^{-1})*
 (E_{\lambda}
 \delta)) = \mu (E_{\lambda} (E_{\lambda} \circ \widetilde{\sigma}))$,
 we conclude that $\lim_{\lambda} \mu(E_{\lambda} \delta^{-1}) =
 \xi^{-1}$, and an argument similar to the one above allows us to draw
 the
 desired conclusion for all negative $n$.
\end{proof}
We may use this to see that $\Xi = \emptyset$ if $(\Delta(A),
 \widetilde{\sigma})$ lacks fixed points.
This is because the restriction of a map $\mu \in \Xi$ to $C_0
 (\Delta(A))$ must be a point evaluation, $\mu_x$ say, by basic Banach
 algebra theory.
If $x \neq \sigma(x)$ there exists an $h \in C_0 (\Delta(A))$ such
 that
 $h(x) = 1$ and $(h\circ \sigma) (x) = 0$.
By Lemma~\ref{apprun} we see that $\mu(h \delta) = \lim_{\lambda} \mu
 (h E_{\lambda} \delta) = \lim_{\lambda} \mu (h) \mu(E_{\lambda}
 \delta))
 = h(x) \xi = \xi$ and likewise $\mu (h \delta^{-1}) = \xi^{-1}$. But
 then $1 = \xi^{-1} \xi = \mu((h \delta^{-1})*(h \delta)) = \mu(h \cdot
 (h
 \circ \sigma)) = h(x) \cdot (h\circ \sigma)(x) = 0$, which is a
 contradiction.

Now for any $x \in \Per^1 (\Delta(A))$ and $\xi \in \mathbb{T}$ there
 is a unique element $\mu \in \Xi$ such that $\mu (f_n \delta^n) = f_n
 (x) \xi^n$ for all $n$ and by the above every element of $\Xi$ must
 be
 of this form for a
 unique $x$ and $\xi$. Thus we have a bijection between $\Xi$ and
 $\Per^1 (\Delta(A)) \times \mathbb{T}$. Denote by $I(x, \xi)$ the
 kernel of
 such $\mu$. This is clearly a modular ideal of
 $\ell_1^{\widehat{\sigma}} (\mathbb{Z}, C_0 (\Delta(A)))$ which is
 maximal and contains $\mathscr{C}$ by multiplicativity and continuity
 of
 elements in
 $\Xi$.

\begin{thm}\label{mamo}
The  modular ideals of $\ell_1^{\widehat{\sigma}} (\mathbb{Z}, C_0
 (\Delta(A)))$ which are maximal and contain the commutator ideal
 $\mathscr{C}$ are precisely the ideals $I(x, \xi)$, where $x \in
 \Per^1
 (\Delta(A))$ and $\xi \in \mathbb{T}$.
\end{thm}
\begin{proof}
One inclusion is clear from the discussion above. For the converse,
 let
 $M$ be such an ideal and note that it is easy to show that a maximal
 ideal containing $\mathscr{C}$ is not properly contained in any
 proper
 left or right ideal. Thus as $\ell_1^{\widehat{\sigma}} (\mathbb{Z},
 C_0
 (\Delta(A)))$ is a spectral algebra, \cite[Theorem 2.4.13]{palmer}
 implies that $\ell_1^{\widehat{\sigma}} (\mathbb{Z}, C_0 (\Delta(A)))
 / M$
 is isomorphic to the complex field. This clearly implies that $M$ is
 the kernel of a non-zero element of $\Xi$.
\end{proof}
\section*{Acknowledgments}
This work was supported by a visitor's grant of the
 \textit{Netherlands
 Organisation for Scientific Research (NWO)}, \textit{The Swedish
 Foundation for International Cooperation in Research and Higher
 Education
 (STINT)}, \textit{Crafoord Foundation, The Royal Swedish Academy of Sciences and The Royal Physiographic Society in Lund}.

\end{document}